\documentclass[12pt]{article}
\usepackage[utf8]{inputenc}
\usepackage[T1]{fontenc}
\usepackage{amssymb}
\usepackage{amsthm}
\usepackage{longtable}
\usepackage{enumitem}

\usepackage{color}
\definecolor{DarkOlive}{rgb}{0.1047,0.2412,0.0064}
\definecolor{FireBrick}{rgb}{0.5812,0.0074,0.0083}
\definecolor{RoyalBlue}{rgb}{0.0236,0.0894,0.6179}
\definecolor{RoyalGreen}{rgb}{0.0236,0.6179,0.0894}
\definecolor{RoyalRed}{rgb}{0.6179,0.0236,0.0894}
\definecolor{LightBlue}{rgb}{0.8544,0.9511,1.0000}
\definecolor{Black}{rgb}{0.0,0.0,0.0}
\definecolor{MidnightBlue}{rgb}{0.0035,0.0020,0.1363}
\definecolor{FireBrick3}{rgb}{0.5812,0.0074,0.0083}
\definecolor{FireBrick4}{rgb}{0.2156,0.0023,0.0035}
\definecolor{Blue2}{rgb}{0.0000,0.0000,0.8644}
\definecolor{Navy}{rgb}{0.0000,0.0000,0.1927}
\definecolor{MediumBlue}{rgb}{0.0000,0.0000,0.6179}
\usepackage[
        a4paper=true,bookmarks=false,pdftitle={Generic Computations in
        Finite Groups of Lie Type},
        colorlinks=true,backref=false,breaklinks=true,linkcolor=MediumBlue,
        citecolor=FireBrick3,filecolor=RoyalRed,
        urlcolor=Blue2,pagecolor=MediumBlue]{hyperref}

\theoremstyle{plain}
\newtheorem{Thm}{Theorem}[section]

\theoremstyle{definition}

\newcommand{\mod}{ \textrm{ mod } }
\newcommand{\Aut}{ \textrm{Aut} }
\newcommand{\F}{\mathbb{F}}
\newcommand{\Z}{\mathbb{Z}}
\newcommand{\R}{\mathbb{R}}

\newcommand{\GL}{\textrm{GL}}

\newcommand{\GF}[1]{\F_{#1}}

\newcommand{\GFp}{\GF{p}}

\begin{document}
\title{Orbits of $Z \circ (2.O_8^+(2).2)$ in Dimension 8}
\author{Frank Lübeck}
\maketitle
\begin{abstract}
Groups  of structure  $2.O_8^+(2)$ have  an irreducible  representation of
degree $8$  which can be  realized over $\Z$  and any prime  field $\GFp$.
This representation  extends to a  group of structure  $2.O_8^+(2).2$. Any
subgroup  $Z \leq  \GFp^{\times}$ acts  by scalar  multiplication on  this
module over $\GFp$.

In this short note we determine for which primes $p > 7$ and which $Z$ the
central products $Z \circ (2.O_8^+(2)$ and $Z \circ (2.O_8^+(2).2)$ have a
regular orbit on the $8$-dimensional $\GFp$-module.

This work was triggered by an omission in the paper~\cite{KP01} by Köhler
and Pahlings,  a paper  which is  used in  various places  in work  on the
$k(GV)$-problem.
\end{abstract}

\section{Introduction}

In~\cite{KP01} Köhler and Pahlings investigated the following problem:

Let $p$ be a prime, $G$ be a  finite group whose order is not divisible by
$p$ and $V$ be a finite  dimensional faithful $\GFp G$-module of dimension
$n$. Furthermore,  assume that $G$  has a quasisimple normal  subgroup $E$
which also acts irreducibly on $V$. Does $G$ have a regular orbit on $V$?

In most cases the  answer to this question is yes, but there  is a list of
exceptions. The  main result of~\cite[Theorem  2.2]{KP01} is the  table of
these exceptions.

For a  fixed quasisimple group  $E$ and  $\GFp E$-module $V$  the possible
groups $G$ (as groups of endomorphisms of $V$) are generated by a subgroup
of  $\Aut(E)$ and  a subgroup  $Z  \leq \GFp^\times$  of scalar  matrices,
see~\cite[Section 3]{KP01}.

From  now  we  consider  the  specific  case  $E  =  2.O_8^+(2)$  and  its
irreducible  representation of  degree  $n  = \dim(V)  =  8$,  hence $p  >
7$.  In this  case  $\Aut(E)$  has the  structure  $2.O_8^+(2).2$ and  the
representation extends in two  ways to this group, see~\cite[p.85]{ATLAS}.
The authors of~\cite{KP01}  determine for which $p$ the groups  of form $Z
\circ E$ have a regular orbit on  $V$ by an elaborate computation with the
table of marks of $O_8^+(2)$. But  they forgot (in statement and proof) to
handle the groups of form $Z \circ  E.2$. In this short note we will close
this gap. We will also recover (with a slight correction) their result for
the groups $Z \circ E$ with an easier argument.

\section{The groups $2.O_8^+(2).2$ and its $8$-dimensional representations
over $\GFp$}

There are  two isomorphism types  of groups with  structure $2.O_8^+(2).2$
which  are  isoclinic, see~\cite[Ch.6,  Sec.7]{ATLAS}.  For  one type  the
(Brauer)-characters of the $8$-dimensional  representations have values in
the rational integers and for the other type the character values generate
$\Z[i]$ ($i^2=-1$). In the latter case the $8$-dimensional representations
can only  be realized over $\GFp$  if the field contains  primitive fourth
roots of unity, that  is if $p \equiv 1 \mod 4$.  If $Z \leq \GFp^\times$,
$|Z|  = 4$,  and $G$  is  one group  of  type $2.O_8^+(2).2$  we find  the
isoclinic group  as subgroup  of index  $2$ in $Z  \circ G$  (exchange the
generators $x$  of $G$  which are  not in the  derived subgroup  $G'$ with
$i\cdot x$).

Now  let $W$  be  the Weyl  group  of  type $E_8$.  It  has the  structure
$2.O_8^+(2).2$,  see~\cite[p.85]{ATLAS}. Since  Weyl groups  have rational
character values it is clear which one of the isoclinic groups this is. We
denote $\tilde W$ the isoclinic group.

Now we can state our result.

\begin{Thm}\label{main}
A  group of  type  $Z \circ  (2.O_8^+(2))$  has no  regular  orbit on  its
$8$-dimensional irreducible $\GFp G$-module if and  only if $p \leq 23$ or
$p=31$ and $|Z| > 2$.

A  group  $G$  of   type  $Z  \circ  W$  has  no   regular  orbit  on  its
$8$-dimensional irreducible $\GFp  G$-modules, if and only if  $p \leq 29$
or $p=31$ and $|Z| > 2$.

In case $p \equiv 1 \mod 4$ a group  $G$ of type $Z \circ \tilde W$ has no
regular orbit on  its $8$-dimensional irreducible $\GFp  G$-modules if and
only if one of the following holds
\begin{itemize}
\item $p < 29$,
\item $p = 29$ and $4 \mid |Z|$,
\item $p = 31$ and $|Z| > 2$.
\end{itemize}
\end{Thm}

\section{The proof}

From one of the irreducible representations of $2.O_8^+(2).2$ of dimension
$8$  we  get   the  other  one  by  tensoring   with  the  $1$-dimensional
representation with kernel $2.O_8^+(2)$. So the action of elements outside
the  derived subgroup  only differs  by scalar  multiplication with  $-1$.
Since the  central element,  which is contained  in the  derived subgroup,
also acts  by the scalar  $-1$, the  orbits on $V$  are the same  for both
module  structures.  Therefore,  it  is  enough to  consider  one  of  the
$8$-dimensional modules.

We now consider the  Weyl group $W$ of type $E_8$. It  can be described as
follows as subgroup of $\GL_8(\Z)$, and this  way we get for any prime $p$
an  $8$-dimensional representation  of  $W$ over  $\GFp$  by reducing  the
matrix entries modulo $p$.

Let $Y = \Z^8$ and $X$ be  the dual $\Z$-lattice. We describe a root datum
of type $E_8$. For  this we take the standard basis vectors  of $Y$ as set
of simple coroots $\alpha^\vee_j$, $1 \leq j \leq 8$, and as corresponding
simple roots  $\alpha_j \in X$ the  rows of the following  matrix (written
with respect to the $\Z$-basis of $X$ which is dual to the simple coroots,
the elements of this basis are also called fundamental weights):
\[ \small
\left(\begin{array}{rrrrrrrr}%
2&0&-1&0&0&0&0&0\\%
0&2&0&-1&0&0&0&0\\%
-1&0&2&-1&0&0&0&0\\%
0&-1&-1&2&-1&0&0&0\\%
0&0&0&-1&2&-1&0&0\\%
0&0&0&0&-1&2&-1&0\\%
0&0&0&0&0&-1&2&-1\\%
0&0&0&0&0&0&-1&2\\%
\end{array}\right)%
\]
For $1 \leq j \leq 8$ we  can use $\alpha_j$ and $\alpha_j^\vee$ to define
the following reflection on $X$:
\[ s_j: X \to X, \quad x \mapsto x - \alpha_j^\vee(x) \cdot \alpha_j.
\]

The group generated  by these reflections $W  = \langle s_j \mid  1 \leq j
\leq 8  \rangle$ is the Weyl  group of type  $E_8$, it is a  Coxeter group
with the $s_j$, $1 \leq j \leq 8$, as set of Coxeter generators. The orbit
$\alpha_1 W$  is called the root  system of type $E_8$,  it contains $240$
roots. The dual  construction on $Y$ yields the  corresponding coroots and
the highest coroot  (the one with is componentwise  $\geq$ the coordinates
of all other coroots) is
\[ \alpha_0^\vee = 
\left(\begin{array}{rrrrrrrr}%
2&3&4&6&5&4&3&2\\%
\end{array}\right)%
.\]

The corresponding root is 
$\alpha_0 = \left(\begin{array}{rrrrrrrr}%
0&0&0&0&0&0&0&1\\%
\end{array}\right)$ 
and this defines a reflection $s_0:X \to X$ as above.

For any prime $p > 7$ we are interested in the orbits of the action of $W$
on $X$ modulo $pX$ (identifying $\GFp = \Z/p\Z$ and $\GFp^8 = X/pX$).

The  subgroup  of  bijections  $X  \to   X$,  generated  by  $W$  and  the
translations by elements of $pX$ is  called the affine Weyl group $W_p$ of
type $E_8$, see~\cite[4.3,  4.8]{Hum90} (we use that the  $\Z$-span of the
roots is all of  $X$, and we scale the translations in  the reference by a
factor $p$).

The action  of $W_p$ on  $X$ can be  $\R$-linearly extended to  the vector
space $X \otimes  \R$. We consider the following subset  of $X \otimes \R$
which is called the  bottom alcove: \[ A_0 := \{ x  \in X \otimes \R\mid\;
\alpha_j^\vee(x) > 0 \mbox{ for } 1\leq j \leq 8, \alpha_0^\vee(x) < p \},
\] and will use the following theorem, see~\cite[4.8]{Hum90}.

\begin{Thm}
The closure $\bar{A}_0$ of $A_0$ is a fundamental domain for the action of
$W_p$  on $X  \otimes \R$.  If $x  \in \bar{A}_0$  then the  stabilizer of
$x$  in  $W$  is generated  by  those  $s_j$,  $1  \leq j  \leq  8$,  with
$\alpha_j^\vee(x) = 0$ together with $s_0$ in case $\alpha_0^\vee(x) = p$.

In  particular,  every  $W_p$-orbit  on   $X  \otimes  \R$  has  a  unique
representative $x \in  \bar{A}_0$ and the orbit is regular  if and only if
$x \in A_0$.
\end{Thm}

Restricting this  to points  $x =  (x_1,\ldots, x_8)  \in X$  we conclude:
There exists a regular  $W$-orbit on $X$ modulo $pX$ if  and only if $A_0$
contains  a  point  with  integer  coordinates.  That  is,  all  $x_j  \in
\Z_{>0}$ (because  $\alpha_j^\vee(x) =  x_i$) for  $1 \leq  j \leq  8$ and
$\alpha_0^\vee(x) < p$.  Since the coordinates of  $\alpha_0^\vee$ are all
positive, such an $x$ exists if and only if
$\rho := \left(\begin{array}{rrrrrrrr}%
1&1&1&1&1&1&1&1\\%
\end{array}\right) \in A_0$ if and only if $\alpha_0^\vee(\rho) = 29 < p$.

It  is an  easy programming  exercise to  enumerate for  moderate $p$  all
integral points $  x \in \bar{A}_0$, and to read  off their stabilizers in
$W$. For all $p  < 29$ these stabilizers all have order  $> 2$. This shows
that also the orbits of the derived  subgroup $W' = 2.O_8^+(2)$ of $W$ are
never regular.

The case  $p=29$. Here all integral  $x \in \bar{A}_0$ with  $x \neq \rho$
have  stabilizer  of  order  $>2$.  But $  \rho  \in  \bar{A}_0$  and  its
stabilizer  is  generated  by  $s_0$  and  is  of  order  $2$.  Since  the
reflections are  a single conjugacy class  of $W$ and generate  $W$ we see
that  $s_0 \notin  W'$. So,  the  orbit of  $\rho$ is  (the only)  regular
$W'$-orbit.

We have shown our Theorem~\ref{main} for $G = W'$ and $G = W$.

\subsection*{Action of scalars}

Now  we consider  the action  of scalars.  We will  see that  we find  all
information we need by considering the orbit of $\rho \in X$ modulo $pX$.

We want  to know for all  primes $p$ all $c  \in \Z$ modulo $p$  such that
there is  an element $w \in  W$ with $\rho w  \equiv c \rho \mod  pX$. The
center of  $W' = 2.O_8^+(2)$ acts  as scalar $-1$. So,  all $W'$-orbits on
$X$ will be closed under multiplication with $-1$.

Fixing $w \in W$ and setting $y = (y_1, \ldots, y_8) := \rho w$ we want to
know for  which $p$ there is  a $c$ such that  $y \equiv c \rho  \mod pX$.
This relation is equivalent  to the condition $\textrm{gcd}(y_1-c, \ldots,
y_8-c) \equiv 0 \mod p$ for some $c \in \Z$. We use
\[\textrm{gcd}(y_1-c, \ldots, y_8-c) = \textrm{gcd}(y_1-c, y_2-y_1, \ldots
y_8-y_1) \mid \textrm{gcd}(y_2-y_1, \ldots, y_8-y_1) \]
and compute the latter expression. 

The possible $p$ are the  prime divisors of $\textrm{gcd}(y_2-y_1, \ldots,
y_8-y_1)$. And  it is  clear that for  each such prime  there is  some $c$
(unique modulo $p$) such that $y_1-c$ is also divisible by $p$.

Using  a computer  and GAP~\cite{GAP}  we computed  the full  $W$-orbit of
$\rho$ and  the decribed $\textrm{gcd}$'s.  For this we used  the explicit
construction of the representation given above. The computation took about
45 minutes on the authors notebook.

The prime divisors  of these numbers are  all $\leq 31$. So, for  $p > 31$
the only $\GFp$- multiple of $\rho$ occuring in the orbit of $\rho$ modulo
$pX$ is $-\rho$, and therefore the orbit of $\rho$ is also regular for any
central product $Z \circ W$ and so for $Z \circ W'$, $Z \leq \GFp^\times$.

The case $p=31$. It is easy to  see that $\rho$ is the only integral point
in $A_0$ for  $p=31$, so there is only one  regular $W$-orbit modulo $31$,
namely  the orbit  of $\rho$.  Therefore it  is not  surprising that  this
$W$-orbit contains all  scalar multiples of $\rho$. Looking  closer at the
set  of $w  \in W$  which yield  multiples of  $\rho$ modulo  $31$ in  our
$\textrm{gcd}$-computations we notice that they are all of even length, so
that already  the $W'$-orbit  of $\rho$ contains  all multiples  of $\rho$
modulo $31$.

This proves the exceptions for $p = 31$ in Theorem~\ref{main}.

The case $p=29$. Our $\textrm{gcd}$-computations show that in the orbit of
$\rho$ modulo $29$ the only multiples  of $\rho$ are $\pm \rho$. The orbit
is not regular, and  $\rho w = \rho$ modulo $29$ for $w=1$  and $w = s_0$.
Now we consider  the isoclinic group $\tilde W$ which  is generated by the
$i s_j$,  $1 \leq j  \leq 8$ where  $i \in \GF{29}$  is of order  $4$. The
orbit of $\rho$  under $\tilde W$ contains powers of  $i$ scalar multiples
of the  vectors in  the orbit  under $W$. The  element $s_0  \in W$  is of
length $57$, so that the corresponding product of generators of $\tilde W$
will map  $\rho$ to $i  \rho$. This shows that  the orbit of  $\rho$ under
$\tilde W$ contains all  four multiples $i^k \rho$, $0 \leq  k \leq 3$. So
$\rho \tilde W$ is twice as long as $\rho W$, hence a regular orbit.

We have  shown all statements  in Theorem~\ref{main} concerning  the cases
with $p=29$. This finishes our proof.

\textbf{Acknowledgement.}  I thank  Melissa  Lee for  making  me aware  of
the  gap concerning  $O_8^+(2)$ in  the paper~\cite{KP01};  this caused  a
corresponding  gap  in early  versions  of  her article~\cite{L20}  (which
hopefully can now be closed).

\bibliographystyle{alpha}
\bibliography{OrbitsWE8}

\end{document}